\documentclass[fleqn]{article}

\usepackage{latexsym}
\usepackage{geometry}
\usepackage{amsmath}
\usepackage{amsthm}
\usepackage{amssymb}
\usepackage{amsfonts}
\usepackage[T1]{fontenc}
\usepackage[rightcaption]{sidecap}

\usepackage{xcolor}
\usepackage{graphicx}

\newcommand{\ods}{\par \vspace{0.2cm} \par}

\newcommand{\dis}{\displaystyle}
\newcommand{\rf}[1]{(\ref{#1})}

\newcommand{\ba}{\begin{array}}
\newcommand{\ea}{\end{array}}

\newcommand{\be}{\begin{equation}}
\newcommand{\ee}{\end{equation}}

\DeclareMathOperator{\sgn}{sgn}

\newcommand{\tta}{\vartheta}

\newtheorem{prop}{Proposition}[section]

\newtheorem{cor}[prop]{Corollary}

\begin{document}

\title{\bf   On geometric interpretation of Euler's substitutions}
\author{{\bf Jan L.\ Cie\'sli\'nski}\thanks{E-mail: \ 
 j.cieslinski@uwb.edu.pl}, \quad  
{\bf  Maciej Jurgielewicz}
\\ {\footnotesize Uniwersytet w Bia\l ymstoku,
Wydzia{\l} Fizyki}
\\ {\footnotesize ul.\ Cio\l kowskiego 1L,  15-245 Bia\l ystok, Poland}
}

\date{ }

\maketitle

\newcommand{\orcidauthorB}{0000-0003-1730-0950} 

\abstract{We consider a classial case of irrational integrals containing a square root of a quadratic polynomial.  It is well known that they can be expressed in terms of elementary functions by one of  three Euler's substitutions. It is less known that the Euler substittutions have a beautiful geometric interpretation.  In the framework of this interpretation one can see that the number 3 is not the most suitable. We show that it is natural to introduce the  fourth Euler substitution. By the way, it is not clear who was the first to attribute these three substitutions to Euler. In his original treatise Leonhard Euler uses two substitutions which are sufficient to cover all cases. 
 }
\ods

\section{Introduction} 

Integrals of rational functions can be expressed in terms of elementary functions.  Therefore a natural method of integration consists in using suitable substitutions and integration by parts to reduce our problem to integration of rational functions. 

In this paper we consider irrational integrals containing the quadratic root of a quadratic polynomial, i.e., integrals of the form  
\be  \label{Rxy}
    \int R (x, y) dx \ ,  
\ee
where $R$ is a rational fuction (a quotient of two polynomials) of $x$ and $y$, and  
\be   \label{y}
   y = \sqrt{a x^2 + b x + c} \ . 
\ee
The subject is, in principle, known. A standard method to deal with such integrals consists in using one of the so called Euler's substitutions \cite{Wiki,Enc}.  However, there are some details which need to be clarified. We will describe in detail a geometric approach to this problem and explain how many Euler substitutions actually exist.  

In fact, to our beest knowledge,  all sources mention exactly three types of substitutions in this context.  It is not clear who was the first to introduce such classification.  Euler substitutions are usually introduced and discussed in Russian sources, see, e.g., \cite{RG,Fikh,Pis}  (Leonhard Euler, although of Swiss origin, lived and worked in Saint Petersburg for many years).  
Surprisingly enough, the three substitutions appeared also in an old textbook by a Harvard professor, \cite{Bye}, without any reference to Euler.  

In our paper we present a clear geometric intepretation of this problem, shortly mentioned in some sources, mainly of Russian origin \cite{Enc,Fich}.  The textbook \cite{Fich} is not translated into English. Another book by the same author, \cite{Fikh}, does not mention this geometric approach in the section on Euler's substitutions.

The main novelty of this paper is  the introduction of the fourth Euler substitution, which is a natural consequence of the geometric approach discussed in our paper.

\ods

\section{Three classical Euler's substitutions}
\label{sec-classical}

The main idea of Euler's substitutions consists in expressing  $ \sqrt{a x^2 + b x + c}$ \ as a linear function of $x$ and a new parameter $t$ in such a  way that the resulting equation is linear with respect to $x$.  In this paper we use the most common numbering of these three substitutions, compare \cite{Wiki,Enc,Fikh,Pis,Bye}. In some sources a different order is used, see \cite{RG,Lu}. 

\subsection{\bf First Euler substitution}
This substitution can be done only in the case $a > 0$:
\be  \label{first}
   \sqrt{a x^2 + b x + c} = \pm x \sqrt{a} + t  \ .   
\ee
Squaring both sides we get:
\[
{a x^2} + b x + c = {a x^2}  \pm 2 x t \sqrt{a} + t^2 \ .
\]
Terms quadratic in $x$ cancel out and the resulting equation is linear in $x$. Computing $x$, we get a rational dependence on $t$:
\be
  x = \frac{t^2 - c}{b \mp 2 t \sqrt{a}} \ .
\ee
Then, from \rf{y} and \rf{first}, we get 
\be
   y = \frac{\mp t^2 \sqrt{a} + t b \mp c \sqrt{a}}{b \mp 2 t \sqrt{a}} \ .
\ee

\subsection{\bf Second Euler substitution}
This substitution can be done only in the case $c > 0$:
\be  \label{second}
   \sqrt{a x^2 + b x + c} = x t \pm \sqrt{c} \ .   
\ee
Squaring both sides we get:
\be
    a x^2 + b x + {c} = x^2 t^2 \pm 2 x t \sqrt{c} + {c} \ .
\ee
The constant $c$ cancels out and dividing both sides by $x$  we again derive an equation linear in $x$. Hence,  similarly as in the previous case, 
\be
 x = \frac{b \mp 2 t \sqrt{c}}{ t^2 - a} \ ,  \qquad  y =  \frac{ b t \mp (t^2 + a) \sqrt{c} }{t^2 - a} \ .
\ee

\subsection{\bf Third Euler substitution}
This substitution can be done only in the case $\Delta > 0$, where 
\be
\Delta \equiv  b^2 -  4 a c 
\ee
 is the discriminant of the quadratic polynomial. Then the polynomial has two distinct real roots $x_1$ and $x_2$, and the third Euler substitution is given by:
\be  \label{third}
   \sqrt{a x^2 + b x + c} = (x - x_1)  t  \ .   
\ee
Squaring both sides we get:
\be
  a (x - x_1)(x-x_2) = (x - x_1)^2  t^2  \qquad \Rightarrow  \qquad   a (x - x_2) = (x - x_1) t^2 \ .
\ee
Computing $x$ from the resulting equation and then using \rf{third} and \rf{y} we obtain 
\be
    x = \frac{t^2 x_1 - a x_2}{t^2 - a} \ , \qquad  y = \frac{ (x_1 - x_2) a t }{t^2 - a} \ ,
\ee
where, of course, 
\be
    x_{1,2} = \frac{ - b \pm \sqrt{b^2 - 4 ac}}{2 a} \ .
\ee

\subsection{\bf Original Euler's approach}

It is interesting that Leonhard Euler himself in his famous monograph used only two of these substitutions, see \cite{Eu}.  He considered two cases:  $\Delta > 0$ and $\Delta < 0$.  
In the first case ($\Delta > 0$) he proposed the substitution \rf{second}, while in the second case ($\Delta < 0$)  he proposed the substitution \rf{first} in a slightly modified form:
\be
    \sqrt{a x^2 + b x + c} = x \sqrt{a} -  t \sqrt{c}  \ .
\ee
Obviously, the case $\Delta  = 0$ is not included because then the quadratic polynomial is a square of the linear function in $x$ and $y$ is linear is $x$ as well. Hence the integrand in \rf{Rxy} is rational in $x$ from the very beginning.

\section{Geometric interpretation}

It is convenient to square both sides of represent \rf{y} obtaining the equation of a quadratic curve
\be  \label{conic}
    y^2 = a x^2 + b x + c
\ee
We will denote this  curve (a conic section) by $Q_{abc}$, i.e., $(x, y) \in Q_{abc}$

\subsection{Elliptic case: $a<0$}

The canonical form of the quadratic polynomial yields:  
\be
      y^2 + |a| \left(  x - \frac{b}{2 |a|} \right)^2 =  c - \frac{b^2}{4 a}
\ee
We can distinguish three cases, depending on the sign of the discriminant $\Delta$:
\be
     \Delta < 0 \quad \Longrightarrow \quad   Q_{abc} = \emptyset \quad \text{(empty set)}
\ee
\be
      \Delta = 0 \quad \Longrightarrow \quad   Q_{abc} = \left\{ \left( \frac{b}{2 |a|}, 0 \right) \right\} \quad \text{(single point)}
\ee
\be
      \Delta > 0 \quad \Longrightarrow \quad   Q_{abc} \text{\ is an ellipse}
\ee
Only in the last case we get a non-degenerated quadratic curve.

\subsection{Parabolic case: $a=0$}

For $a=0$ (and $b\neq 0$)  the conic  $Q_{abc}$ is a parabola with the symmetry axis $y=0$. 

\subsection{Hyperbolic case:  $a>0$}
The canonical form of the quadratic polynomial yields:  
\be
      y^2 - a \left(  x + \frac{b}{2 a} \right)^2 =  c - \frac{b^2}{4 a}
\ee
We can distinguish three cases, depending on the sign of the discriminant $\Delta$:
\be
     \Delta < 0 \quad \Longrightarrow \quad   Q_{abc} \text{\ is a hyperbola with vertices at the line $x = - \frac{b}{2 a}$} 
\ee
\be
      \Delta = 0 \quad \Longrightarrow \quad   Q_{abc} \text{\ is a pair of intersection lines} 
\ee
\be
      \Delta > 0 \quad \Longrightarrow \quad   Q_{abc} \text{\ is a hyperbola with vertices at $x$ axis}
\ee
Therefore, for $\Delta \neq 0$ we get a non-degenerated quadratic curve.

\subsection{Rational parameterization -- standard approach} 

The key idea leading to a rational parameterization consists in fixing an arbitrary point $P_0 = (x_0, y_0)$ on the conic $Q_{abc}$ and assigning to any other point $P=(x, y)$ of this conic the line $P_0 P$.  Taking as a parameter $t$ the slope of this line  we obtain a rational parameterization of the  conic $Q_{abc}$ \cite{Enc,Fich}. Thus we have the system of three equations:
\be  \ba{l}
  y^2 =  a x^2 + b x + c \ , \\[1ex]
 y_0^2 = a x_0^2 + b x_0 + c , \\[1ex]
  y - y_0 = t (x - x_0) \ . 
\ea \ee
The points $(x,y)$ and $(x_0, y_0)$ belong to the conic $Q_{abc}$ and $t$ is the slope of the straight line passing through $(x,y)$ and $(x_0, y_0)$.  Subtracting the second equation from the first one we get:
\be  \ba{l}
    (y-y_0)(y+y_0) = (x-x_0) ( a (x+x_0) + b ) \ , \\[1ex]
   y_0^2 = a x_0^2 + b x_0 + c , \\[1ex]
 y - y_0 = t (x - x_0) \ , 
\ea \ee
Substituting the last equation into the first one we obtain: 
\be  \ba{l}
   \left( t (y+y_0) -  a (x+x_0) -  b \right) (x - x_0) = 0 \ , \\[1ex]
   y_0^2 = a x_0^2 + b x_0 + c , \\[1ex]
 y - y_0 = t (x - x_0) \ .
\ea \ee
Assuming $x \neq x_0$, we get
\be  \ba{l}
    t (y+y_0) =  a (x+x_0) +  b  \ , \\[1ex]
   y_0^2 = a x_0^2 + b x_0 + c , \\[1ex]
 y - y_0 = t (x - x_0) \ .
\ea \ee
Now, the first and the last equation form a system of two linear equations for two variables $x, y$, which can be solved in the standard way. As a result we obtain:
\be \label{xy gen}  \ba{l}  \dis
x = \frac{x_0 t^2 - 2 y_0 t + a x_0 + b}{t^2 - a} \ , \\[2ex]\dis
y = \frac{ - y_0 t^2 + (2 x_0 + b) t - a y_0}{t^2 - a} \ ,
\ea \ee
which means that we expressed $x$ and $y$ as rational functions of the parameter $t$. 

\begin{cor}
There are infinitely many Euler-like substitutions. Each of them is determined by the choice of $x_0$, provided that $a x_0^2 + b x_0 + c \geqslant 0$. Then the point $P_0 \equiv (x_0, y_0)$ is given by:
\be
     P_0 =  \left(  x_0, \pm \sqrt{a x_0^2 + b x_0 + c} \right)
\ee
\end{cor}

\begin{figure}
\includegraphics{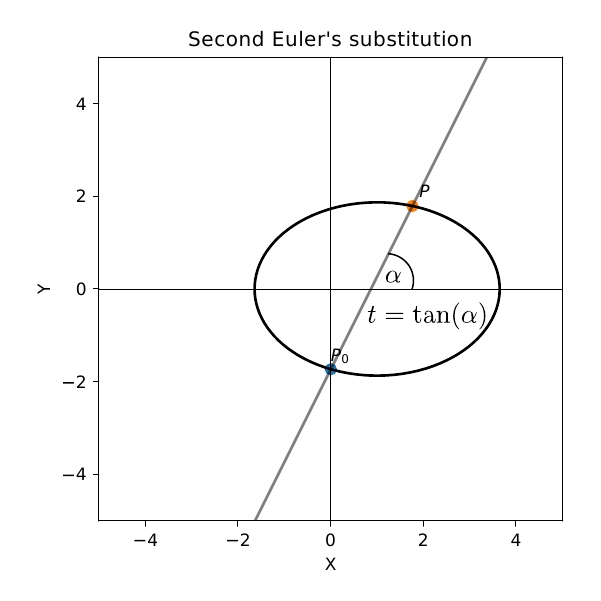}  

\caption{Geometric interpretation of the second Euler substitution in the case $a<0$ and $c>0$. The point $P$  is parameterized by the slope $t$ of the line $P_0 P$. }  
\label{2sub-e}
\end{figure}

\begin{figure}

\includegraphics{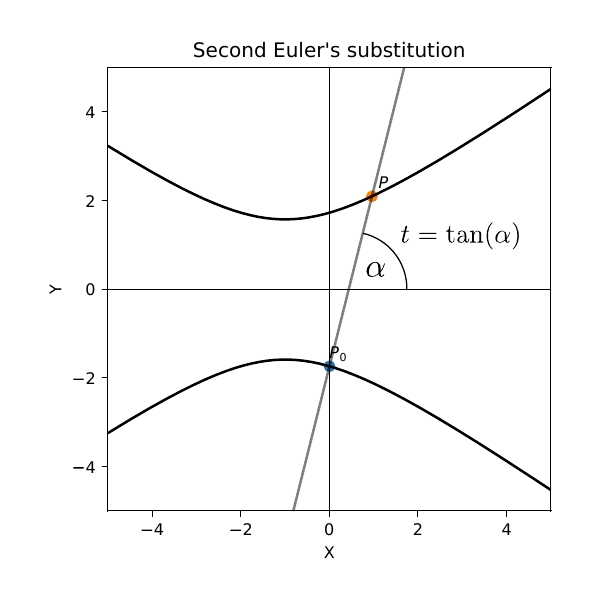}  

\caption{Geometric interpretation of the second Euler substitution in the case $a>0$ and $c>0$. The point $P$  is parameterized by the slope $t$ of the line $P_0 P$. }  
\label{2sub-h}
\end{figure}

\begin{figure}
\includegraphics{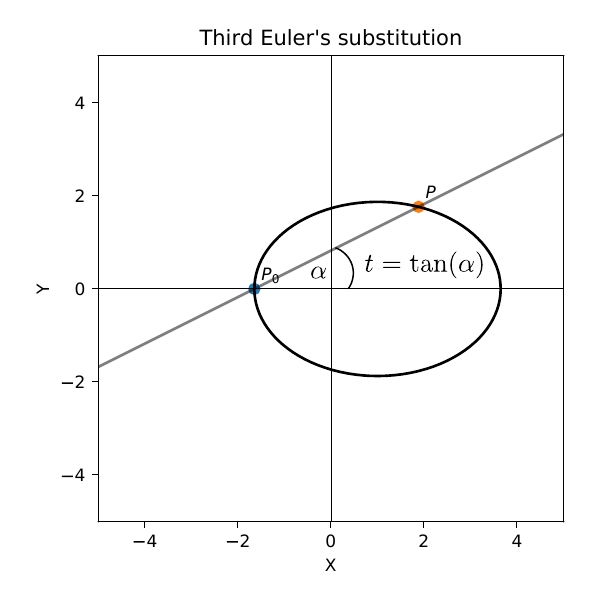}  

\caption{Geometric interpretation of the third Euler substitution in the case $a<0$ and $\Delta>0$. The point $P$  is parameterized by the slope $t$ of the line $P_0 P$. }  
\label{3sub-e}
\end{figure}

\begin{figure}

\includegraphics{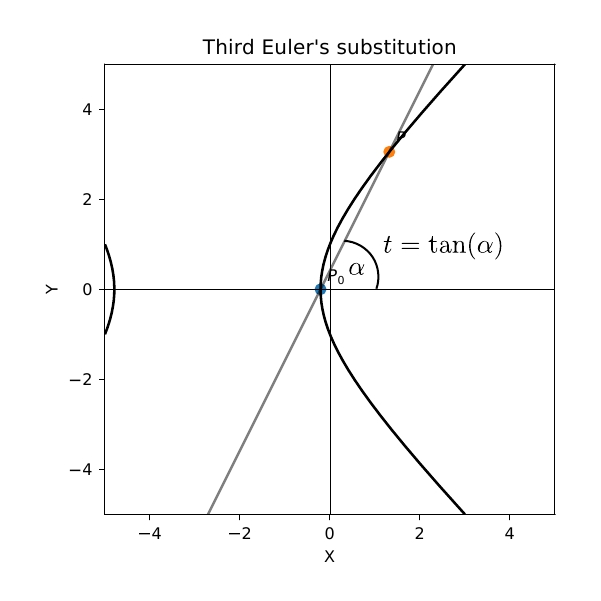}  

\caption{Geometric interpretation of the third Euler substitution in the case $a>0$. The point $P$ is parameterized by the slope $t$ of the line $P_0 P$.  }
\label{3sub-h}

\end{figure}

In particular, the second Euler substitution corresponds to $x_0 = 0$ (provided that the graph of the quadric $Q_{abc}$ intersects the axis $y$), see Fig.~\ref{2sub-e} and Fig.~\ref{2sub-h}. The third Euler substitution corresponds to $x_0$ being a root of the polynomial $a x^2 + b x + c$ (provided that the graph of $Q_{abc}$ intersects the axis $x$), see  Fig.~\ref{3sub-e} and Fig.~\ref{3sub-h}.

\begin{figure}
\includegraphics{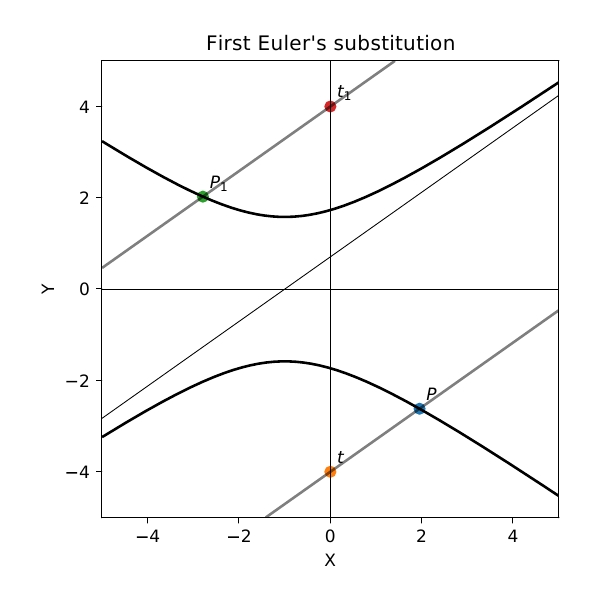}  

\caption{Geometric interpretation of the first Euler substitution. The points $P$ and $P_1$ are parameterized by intersections $t$ and $t_1$, respectively, of  the $y$-axis with the line parallel to one of the asymptotes of the hyperbola $y^2 = a x^2 + b x + c$.  }    \label{1sub-h}

\end{figure}

The first Euler substitution apparently does not fit this picture. However, its geometric interpretation is even simpler and more evident.  The formula \rf{first}  describes the family of lines parallel to asymptotes of the corresponding hyperbola, see Fig.~\ref{1sub-h}.  We may treat it as a special case of \rf{xy gen} when the point $(x_0, y_0)$ lies at infinity. Note that points $(x_0, \pm x_0 \sqrt{a})$  belong to the conic \rf{conic} in the limit for $x_0 \rightarrow \infty$.    


\begin{figure}
\vspace{-3ex}

\includegraphics{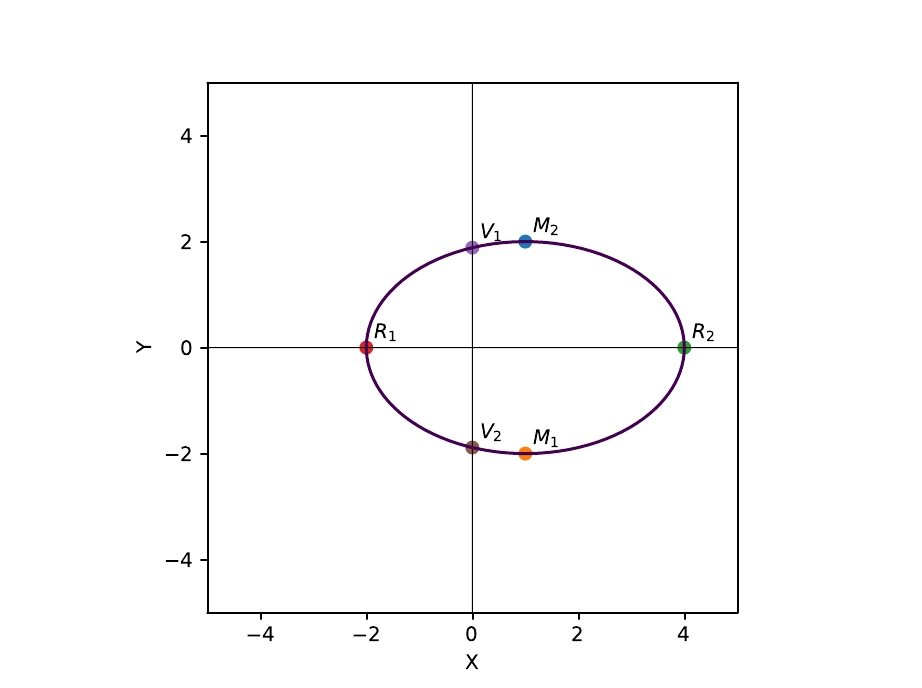}  
\caption{Characteristic points on the graph of an ellipse:  intersections with the coordinate axes (provided that they exist) and extremes (minimum $M_1$ and maximum $M_2$).}  
\label{elipsa}

\vspace{1ex}

\includegraphics{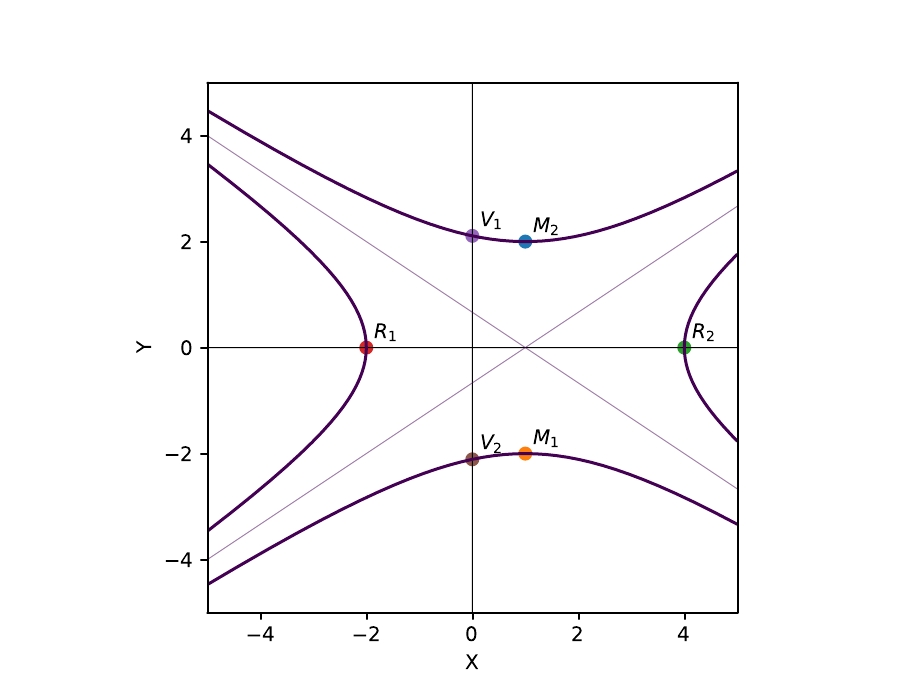}  

\caption{Characteristic points on the graphs of hyperbolas (two hyperbolas with the same $|a|$ are presented):   intersections with the coordinate axes ($V_1$, $V_2$, $R_1$, $R_2$) and extremes ($M_1$, $M_2$).  }     \label{hip}

\end{figure}

\begin{figure}[t]

\includegraphics{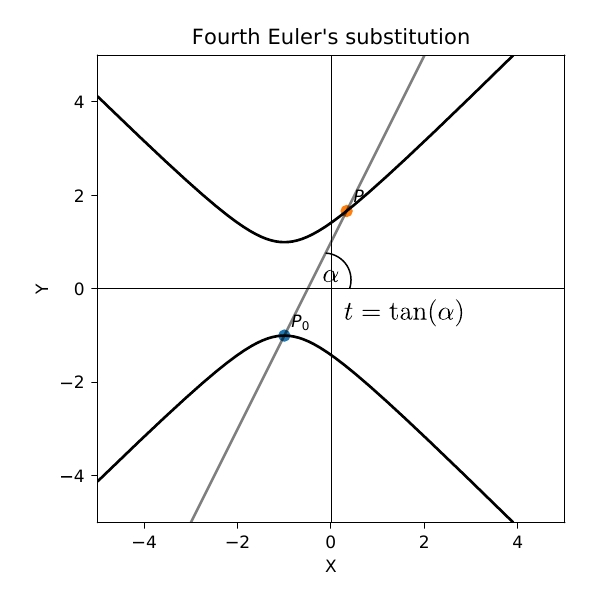}  

\caption{Geometric interpretation of the fourth Euler substitution in the case $a>0$. The point $P$ is parameterized by the slope $t$ of the line $P_0 P$, where $P_0 = M_1$.  }
\label{4sub-h}
\ods

\end{figure}

\begin{figure}

\includegraphics{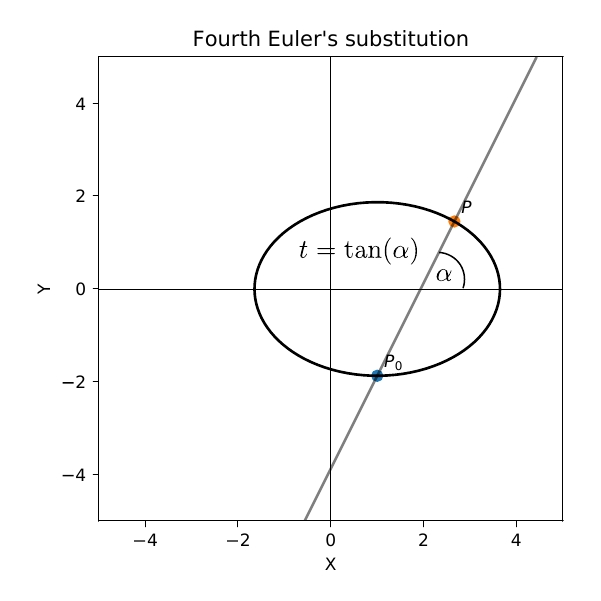}  

\caption{Geometric interpretation of the fourth Euler substitution in the case $a<0$. The point $P$ is parameterized by the slope $t$ of the line $P_0 P$, where $P_0 = M_1$.  }
\label{4sub-e}
\ods

\end{figure}

\ods


\section{New insights from the geometric interpretation}

The description given in the previous section is more or less known (see, e.g., \cite{Enc,Fich}), although we are not aware about any  reference containing all these details. We are going to derive from this geometric picture more quite interesting consequnces.

First of all, we identify characteristic points on the graph of a quadratic curve which can be chosen as $P_0$ in the most natural way: vertices ($M_1$, $M_2$, $R_1$, $R_2$)  and intersections with coordinate axes ($R_1$, $R_2$, $V_1$, $V_2$),  see Fig.~\ref{elipsa} and Fig.~\ref{hip}. 

In particular,  in the case of the second Euler substitution    $P_0 = V_2$  (see Fig.~\ref{2sub-e} and Fig.~\ref{2sub-h}) or $P_0=V_1$, while in the case of the third Euler substitution $P_0 = R_1$ (see Fig.~\ref{3sub-e}) or $P_0=R_2$ (see Fig.~\ref{3sub-h}).  The first Euler substitution is related to  $P_0$ at infinity.

\subsection{\bf Fourth Euler's substitution}

The geometric approach presented above includes all three classical Euler's substitutions, but it is still missing vertices $M_1$ and $M_2$.   Therefore, it is natural to introduce another (fourth) Euler's substitution, geometrically related to missing vertices:  $P_0 = M_1$  (see Fig.~\ref{4sub-e}  and Fig~\ref{4sub-h}) or $P_0 = M_2$. 

The algebraic description of the fourth Euler substitution is based on the canonical form of the quadratic polynomial:
\be  \label{ypq}
y = \sqrt{ a (x - p)^2 + q} 
\ee
where  
\be
  p = - \frac{b}{2 a} \ , \qquad  q = - \frac{\Delta}{4 a} \ . 
\ee
The fourth Euler substitution is defined by:
\be  \label{yq}
  y =  \sqrt{q} + (x-p) t \ .
\ee
Squaring both sides we get:
\be 
    a (x-p)^2 + {q} =  {q} + 2 (x - p) t \sqrt{q} + (x-p)^2 t^2  .
\ee
The constant $q$ cancels out and dividing both sides by $x-p$, we obtain 
\be
    a (x-p) = 2 t \sqrt{q} +(x-p) t^2 , 
\ee
which is linear in $x$. Hence
\be 
x-p = \dis \frac{2 t \sqrt{q}}{a - t^2}   , 
\ee
and using \rf{yq} we get 
\be
y = \frac{a + t^2}{a - t^2} \sqrt{q}  \ .
\ee
Thus we have a rational dependence of $x$ and $y$ on the parameter $t$. Moreover,
\be
\frac{dx}{dt} =  \frac{2 (a + t^2) \sqrt{q}}{(a - t^2)^2}    \ ,        
\ee
and we can easily transform the irrational integral  function  \rf{Rxy} into an integral function rational with respect to $t$. 

\subsection{Rational parameterization -- other parameters}

Geometric approach suggests also some modifications or new variants of the existing rational parameterizations.  
Here we confine ourselves to one example. Introducing a new parameter $\tau$
\be
   \tau = 2 t \sqrt{a} - b\ ,
\ee
we obtain the following simplification of the first Euler substitution:
\be 
x = - \frac{1}{4 a} \left(  \tau + \frac{\Delta}{\tau}  + 2 \right)  \ , \qquad  
y = \frac{1}{4\sqrt{a}}  \left(  \tau - \frac{\Delta}{\tau}  \right) \ .
\ee
Geometrically, the parameter $\tau = 0$  corresponds to the line passing through the point  $(p,0)$ and this is one of two asymptotes (that is why $x \rightarrow \infty$ and $y \rightarrow \infty$ for $\tau \rightarrow 0$).

\section{Euler's substitutions {\it versus} trigonometric substitutions}

Another popular method for computing irrational integrals  \rf{Rxy}  consists in making a suitable trigonometric or hyperbolic substitution.  We use the canonical form of the quadratic curve (compare \rf{ypq}):
\be  \label{y2pq}
 y^2 = a (x -p)^2 + q \ .
\ee
Assuming $q \neq 0$ (otherwise $y$ depends linearly on $x$) we introduce new variables $\xi, \eta$ as folows:
\be
   \eta = \frac{y}{\sqrt{|q|}} \ , \qquad     \xi = \frac{(x - p) \sqrt{|a|}}{\sqrt{|q|} } \ .
 \ee
Then  \rf{y2pq}  becomes
\be
     \eta^2 = (\sgn a) \ \xi^2 + \sgn q \ ,
\ee
because $a/|a| = \sgn a$, etc. 

Thus we have three separate cases (in the fourth case --both signs negative-- there are no real solutions), where trigonometric or hyperbolic substitutions are well known: 
\be \ba{l}
\eta = \sqrt{\xi^2 -1}    \quad  \Longrightarrow \quad   \xi = \cosh \tta  \ , \ \eta = \sinh \tta \ , \\[1ex]
\eta = \sqrt{1 - \xi^2}   \quad  \Longrightarrow \quad   \xi = \cos \tta \ , \ \eta = \sin \tta \ ,  \\[1ex]
\eta = \sqrt{\xi^2 +1}   \quad  \Longrightarrow \quad   \xi = \sinh \tta \ , \  \eta = \cosh \tta \ .
\ea \ee

Is it better than Euler's substitutions? This is a matter of taste.  Perhaps it is more easy to memorize, however, one has to remember that integrals of trigonometric or hyperbolic functions have to be converted into integrals of rational functions by another substitution:
\be
     t = \tan \frac{\theta}{2} \quad \text{or} \quad  t = \tanh\frac{\theta}{2} \ .
\ee

\section{Conclusions}

We presented and discussed a geometric approach to Euler substitutions. One  consequence of this thorough discussion was the introduction of the fourth Euler substitution in addition to three traditionally mentioned Euler substitutions.  In fact, we can say about infinite number (one parameter family) of Euler-like substitutions.  They can be further modified or simplified by suitable linear or fractional linear transformations.

\vspace{6pt}

\end{document}